\documentclass{amsart}
\usepackage{amsfonts}
\usepackage{amsmath}
\usepackage{amssymb}
\usepackage{amsthm}

\newtheorem{theorem}{Theorem}
\newtheorem{lemma}{Lemma}
\newtheorem*{prob}{Problem}

\begin{document}
\author{Mark Pankov}
\title{Order preserving transformations of the Hilbert grassmannian}
\address{Department of Mathematics and Information Technology,
University of Warmia and Mazury,
{\. Z}olnierska 14A, 10-561 Olsztyn, Poland}
\email{pankov@matman.uwm.edu.pl}
\subjclass[2000]{}

\begin{abstract}
Let $H$ be a separable real Hilbert space.
Denote by ${\mathcal G}_{\infty}(H)$ the Grassmannian consisting
of closed subspaces with infinite dimension and codimension.
This Grassmannian is partially ordered by the inclusion relation.
We show that every order preserving transformation of ${\mathcal G}_{\infty}(H)$
can be extended to an automorphism of the lattice of closed subspaces of $H$.
It follows from Mackey's result \cite{Mackey} that automorphisms of this lattice
are induced by invertible bounded linear operators.
\end{abstract}

\maketitle

Dedicated to G. W. Mackey (1916 -- 2006)

\section{Introduction}
Let $V$ be a left vector space over a division ring $R$.
If $\dim V=n$ is finite then we write ${\mathcal G}_{k}(V)$
for the Grassmannian consisting of $k$-dimensional subspaces of $V$,
two elements of ${\mathcal G}_{k}(V)$ are called {\it adjacent} if
their intersection is $(k-1)$-dimensional.
Chow's theorem \cite{Chow} states that
every adjacency preserving transformation of
${\mathcal G}_{k}(V)$ ($1<k<n-1$) is induced by a
semi-linear isomorphism of $V$ to itself or to the dual space $V^{*}$,
and the second possibility can be realized only for the case when $n=2k$.
Recall that a mapping $l:V\to V$ is semi-linear
if it is additive
and there exists an automorphism $\sigma:R\to R$ such that
$$l(ax)=\sigma(a)l(x)$$
for all $x\in V$ and $a\in R$.
Some interesting results closely connected with Chow's theorem
can be found in \cite{Havlicek,Huang,Kreuzer}.

Now suppose that $V$ is infinite-dimensional.
In this case we define Grassmannians as the orbits of
the action of the linear group on the set of proper subspaces of $V$.
Let us consider a Grassmannian ${\mathcal G}$
consisting of subspaces with infinite dimension and codimension.
Two subspaces $S,U\in {\mathcal G}$ are said to be {\it adjacent}
if
$$\dim(S/(S\cap U))=\dim(U/(S\cap U))=1.$$
The Grassmann graph
(the graph whose vertex set is ${\mathcal G}$
and whose edges are pairs of adjacent subspaces)
is not connected and there are adjacency preserving
transformations of ${\mathcal G}$ which are not induced by semi-linear mappings
(see Example 4.3 in \cite{BlunckHavlicek}).
On the other hand, ${\mathcal G}$ is partially ordered by the inclusion
relation and every order preserving transformation of ${\mathcal G}$
preserves the adjacency relation,
however the adjacency preserving transformation constructed in \cite{BlunckHavlicek}
is not order preserving.

In the present paper we consider a separable real Hilbert space
and the Grassmannian consisting of its closed subspaces with infinite dimension
and codimension.
We show that every order preserving transformation of this Grassmannian
can be extended to an automorphism of the lattice of closed subspaces of the Hilbert space.
By Mackey's result \cite{Mackey},
automorphisms of this lattice are induced by invertible bounded linear operators.

\section{Result}
Let $H$ be a separable real Hilbert space.
Denote by ${\mathcal G}(H)$ the lattice of closed subspaces of $H$.
The group ${\rm GL}(H)$ (the group of invertible bounded linear operators)
acts on ${\mathcal G}(H)$, the orbits of this action will be called
{\it Grassmannians}.
There are the following three types of Grassmannians:
\begin{enumerate}
\item[$\bullet$]
${\mathcal G}_{k}(H)$ consisting of $k$-dimensional subspaces,
\item[$\bullet$]
${\mathcal G}^{k}(H)$ consisting of closed subspaces with
codimension $k$,
\item[$\bullet$]
${\mathcal G}_{\infty}(H)$ consisting of closed subspaces
with infinite dimension and codimension.
\end{enumerate}
The latter Grassmannian is partially ordered by the inclusion relation.

Every invertible bounded linear operator induces an automorphism of the lattice
${\mathcal G}(H)$.

\begin{theorem}[\cite{Mackey}, p. 246]
Every automorphism of the lattice ${\mathcal G}(H)$ is induced
by an invertible bounded linear operator.
\end{theorem}

Let $f$ be an automorphism of ${\mathcal G}(H)$.
The restriction of $f$ to ${\mathcal G}_{1}(H)$
is a collineation of the projective space associated with $H$ to itself
(points of this projective space are $1$-dimensional subspaces and lines are defined by
$2$-dimensional subspaces).
This collineation is induced by an invertible linear operator $A:H\to H$
(the Fundamental Theorem of Projective Geometry \cite{Artin, Baer}),
and it is not difficult to show that
$$f(S)=A(S)$$
for every $S\in {\mathcal G}(H)$.
Theorem 1 is  a consequence of the following lemma.

\begin{lemma}[Lemma B in \cite{Mackey}]
Every invertible linear operator preserving ${\mathcal G}^{1}(H)$ is bounded.
\end{lemma}

The results given above were established in \cite{Mackey}
for normed vector spaces over ${\mathbb R}$.

In this paper the following statement will be proved.

\begin{theorem}
Let $f:{\mathcal G}_{\infty}(H)\to {\mathcal G}_{\infty}(H)$
be an order preserving bijection:
$$S\subset U\;\Longleftrightarrow\; f(S)\subset f(U).$$
Then $f$ can be extended to an automorphism of ${\mathcal G}(H)$.
\end{theorem}

\section{Proof of Theorem 2}
First we give a trivial remark concerning the automorphisms of ${\mathcal G}(H)$
induced by invertible bounded linear operators.
Let $A\in {\rm GL}(H)$.
Recall that the {\it adjoint} operator $A^{*}:H\to H$ is defined by the formula
$$(A^{*}x, y)=(x,Ay)$$
for all $x,y\in H$.
An easy verification shows that the mapping
$$S\to (A(S^{\perp}))^{\perp}$$
($S^{\perp}$ is the orthogonal complement to $S$)
is the automorphism of ${\mathcal G}(H)$ induced by the operator $(A^{*})^{-1}$.

\begin{lemma}
Let $S\in {\mathcal G}_{\infty}(H)$
and ${\mathcal X}$ be the set of all elements of ${\mathcal G}_{\infty}(H)$
contained in $S$.
There exists an invertible bounded linear operator $A:S\to f(S)$ such that
$$f(U)=A(U)$$
for all $U\in {\mathcal X}$.
\end{lemma}

\begin{proof}
We restrict ourself to the case when $f(S)=S$,
since in the general case we can take $C\in {\rm GL}(H)$
which sends $f(S)$ to $S$ and consider the transformation $U\to C(f(U))$.

Let us consider $S$ as a Hilbert space and denote by
${\mathcal Y}$ the set of its closed subspaces with infinite codimension.
We have $U\in{\mathcal Y}$ if and only if $U^{\perp}$ belongs to ${\mathcal X}$
(throughout the proof we will write $U^{\perp}$ for the orthogonal complement in $S$).
Since $f$ is order preserving, the same holds for
the transformation $g:{\mathcal Y}\to {\mathcal Y}$ sending
$U$ to $(f(U^{\perp}))^{\perp}$;
and the equality
$$f({\mathcal G}^{i}(S))={\mathcal G}^{i}(S)$$
implies that
$$g({\mathcal G}_{i}(S))={\mathcal G}_{i}(S),$$
in particular the restriction of $g$ to ${\mathcal G}_{1}(S)$
is a collineation of the projective space associated with $S$ to itself.
This collineation is induced by an invertible linear operator $B:S\to S$
(the Fundamental Theorem of Projective Geometry) and
an easy verification shows that
$$g(U)=B(U)$$
for all $U\in {\mathcal Y}$.
For every infinite-dimensional subspace $U\in {\mathcal Y}$  we have
$$g({\mathcal G}^{1}(U))={\mathcal G}^{1}(g(U))$$
($g$ is order preserving) and  Lemma 1 guarantees that
the restriction of $B$ to every element of ${\mathcal Y}$
is bounded.
Since there are $M,N\in {\mathcal Y}$ such that $S=M\dotplus N$,
the operator $B$ is bounded.
By the remark given above, the operator $(B^{*})^{-1}$ is as required.
\end{proof}

\begin{lemma}
Let $S$ and $U$ be elements of ${\mathcal G}_{\infty}(H)$
such that $S\dotplus U$ belongs to ${\mathcal G}_{\infty}(H)$.
Then
$$\dim(S\cap U)=\dim(f(S)\cap f(U)).$$
\end{lemma}

\begin{proof}
We apply Lemma 2 to $S\dotplus U$ and get the claim.
\end{proof}

If $S$ is a finite-dimensional subspace of $H$ then we write $[S]$ for
the set of all elements of ${\mathcal G}_{\infty}(H)$ containing $S$.

\begin{lemma}
For every $S\in {\mathcal G}_{k}(H)$ there exists $S'\in {\mathcal G}_{k}(H)$
such that
$$f([S])=[S'].$$
\end{lemma}

\begin{proof}
We choose $M,N\in [S]$ such that $M\cap N=S$ and define
$$S':=f(M)\cap f(N).$$
By Lemma 3, $S'$ belongs to ${\mathcal G}_{k}(H)$.

Let $U$ be an element of $[S]$.
First we consider the case when $U$ has an infinite-dimensional
intersection with $M$ or $N$.
Suppose that $U\cap N$ is infinite-dimensional.
Then $U\cap N$ belongs to $[S]$ and its $f$-image is contained in $f(N)$.
Hence
$$f(U\cap N)\cap f(M) \subset S'.$$
On the other hand, $U\cap N$ intersects $M$ precisely by $S$.
Lemma 3 guarantees that $f(U\cap N)\cap f(M)$ is $k$-dimensional,
thus this subspace coincides with $S'$.
This implies that
$$S'\subset f(U\cap N)\subset f(U)$$
and $f(U)$ belongs to $[S']$.

Now suppose that $U\cap M$ and $U\cap N$ both are finite-dimensional.
Then $U$ contains $T\in {\mathcal G}_{\infty}(H)$ such that
$$T\cap M=T\cap N=S.$$
By Lemma 3, the subspaces
$$f(T)\cap f(M)\;\mbox{ and }\;f(T)\cap f(N)$$
are $k$-dimensional.
If one of these subspaces
does not coincide  with $S'$
then they are different
(it follows from the equality $f(M)\cap f(N)=S'$).
In this case $f(T)$ contains $T'\in {\mathcal G}_{\infty}(H)$ such that
$$\dim(f(M)\cap T')=k\;\mbox{ and }\;\dim(f(N)\cap T')<k.$$
It follows from Lemma 3 that $f^{-1}(T')$ intersects $M$ by a $k$-dimensional subspace;
since $f^{-1}(T')$ is contained in $T$ and $T\cap M=S$,
this subspace coincides with $S$.
The inclusion $S\subset f^{-1}(T')\subset T$ and the equality  $T\cap N=S$ show that
$$N\cap f^{-1}(T')=S$$
which contradicts the latter inequality.
Therefore
$$f(T)\cap f(M)=f(T)\cap f(N)=S';$$
in particular $f(T)$ is an element of $[S']$,
and the same holds for $f(U)$ (since $f(T)$ is contained in $f(U)$).

We have established that $f([S])\subset [S']$.
The proof of the inverse inclusion is similar.
\end{proof}

If $S\in {\mathcal G}(H)$ is a subspace with finite codimension then we write $[S]$
for the set of all elements of ${\mathcal G}_{\infty}(H)$ contained in $S$.

\begin{lemma}
For every $S\in {\mathcal G}^{k}(H)$ there exists $S'\in {\mathcal G}^{k}(H)$
such that
$$f([S])=[S'].$$
\end{lemma}

\begin{proof}
Let us consider the order preserving transformation
$g:{\mathcal G}_{\infty}(H)\to {\mathcal G}_{\infty}(H)$
sending $U$ to $(f(U^{\perp}))^{\perp}$.
Lemma 4 implies the existence of $T\in {\mathcal G}_{k}(H)$
such that
$$g([S^{\perp}])=[T].$$
The subspace $S'=T^{\perp}$ is as required.
\end{proof}

If $S\in {\mathcal G}(H)$ is a subspace of finite dimension or codimension then
there exists a subspace $S'$ belonging to the Grassmannian containing $S$
and such that
$$f([S])=[S']$$
(Lemmas 4 and 5).
We define $f(S):=S'$.
This gives the required extension.

\section{Problem}

Let $H$ be a separable complex Hilbert space.
A mapping $A:H\to H$ is said to be a {\it semi-linear operator}
if it is additive and there exists an automorphism $\sigma:{\mathbb C}\to {\mathbb C}$
such that
$$Aax=\sigma(a)Ax$$
for all $x\in H$ and all $a\in {\mathbb C}$;
if $A\not \equiv 0$ then there is only one automorphism satisfying
this condition, and we get a usual linear operator if this automorphism is identical.

The field ${\mathbb R}$ has not non-trivial automorphisms,
thus in the real case any semi-linear operator is linear.
The complex case is more complicated:
the conjugate transformation is unique non-trivial continuous automorphism,
but there exist non-continuous automorphisms of ${\mathbb C}$.

If a semi-linear operator $A:H\to H$ ($A\not \equiv 0$)
is bounded then the associated automorphism $\sigma:{\mathbb C}\to {\mathbb C}$
is continuous, hence $\sigma$ is  identical ($A$ is linear) or
it is the conjugate transformation ($Aax=\overline{a}Ax$).

\begin{prob}
Show that every semi-linear operator preserving ${\mathcal G}^{1}(H)$
is bounded.
\end{prob}

The methods used to prove Lemma 1 (Lemma B in \cite{Mackey})
can not be exploited in the complex case, it is connected with
the existence of non-continuous automorphisms of ${\mathbb C}$.

\end{document}